\documentclass[12pt]{amsart}

\usepackage[left=2.8cm, right=2.8cm]{geometry}
\usepackage{stackrel}
\usepackage{amssymb}
\usepackage{amsmath}
\usepackage{amsthm}
\usepackage[pdftex]{graphicx}
\usepackage{mathtools}
\usepackage{pdfsync}
\usepackage{graphicx}
\usepackage{color}
\usepackage{tikz-cd}
\usepackage{tikz}
\usepackage{marginnote}
\usetikzlibrary{positioning}
\usetikzlibrary{decorations.text}
\usetikzlibrary{decorations.pathmorphing}
\usetikzlibrary{matrix}

\usepackage[all]{xy}

\newtheorem{theorem}{Theorem}[section]
\newtheorem{proposition}[theorem]{Proposition}
\newtheorem{lemma}[theorem]{Lemma}
\newtheorem{corollary}[theorem]{Corollary}

\theoremstyle{definition}

\newtheorem{definition}[theorem]{Definition}
\newtheorem{example}[theorem]{Example}

\newtheorem*{remark}{Remark}

\numberwithin{equation}{section}
\DeclarePairedDelimiterX\Set[2]{\lbrace}{\rbrace}
 { #1 \,\delimsize|\,\mathopen{} #2 }

\def\({\left(}
\def\){\right)}

 \newcommand{\R}{\mathbb{R}}
  
 \newcommand{\Z}{\mathbb{Z}}

\title[]{Embeddings of Relative Dimension Groups for Zero-Dimensional Dynamical Systems}

\author{Maryam Hosseini}
\address{Department of Mathematics, School of Science and Technology, City St. Georges University of London, UK }
\email{maryam.hosseini@city.ac.uk}

\begin{document}
\maketitle
%%%%%%%%%%%%%%%%%%%%%%%%%%%%%%%%%%%%%%%%%%%%%%%%%%
\begin{abstract}
Let $\pi:(X,T)\rightarrow(Y,S)$ be a factor map between zero-dimensional dynamical systems, and let $X_0$ and $Y_0$ be quasi-sections with $\pi(X_0)=Y_0$. We study the homomorphism $\pi^*:K^0(Y|Y_0,S)\rightarrow K^0(X|X_0,T)$, $[f]\mapsto[f\circ\pi]$, between the corresponding relative dimension groups, and identify it with the homomorphism induced by the premorphism between Bratteli--Vershik realizations associated with $\pi$. Generalizing a result of Glasner and Weiss for minimal Cantor systems, we show that if $Y_0$ is quasi-simple, then $\pi^*$ is an injective unital order-preserving homomorphism, and that if $X_0=\pi^{-1}(Y_0)$, then $\pi^*$ is a unital order embedding; the latter condition is characterized in terms of the premorphism. Examples show that quasi-simplicity cannot be dropped, and that without the condition $X_0=\pi^{-1}(Y_0)$ the map $\pi^*$ need not be an order embedding.
\end{abstract}
%\tableofcontents{}
%%%%%%%%%%%%%%%%%%%%%%%%%%%%%%%%%%%%%%%%%%%%%%%%%
\section{Introduction}
The interplay between operator algebras and dynamical systems has created a robust framework for the study of dynamical systems. For minimal Cantor systems, Giordano, Putnam and Skau showed that $K$-theory classifies such systems up to orbit equivalence (i.e., there is a homeomorphism from $X$ onto $Y$ mapping each $T$-orbit onto an $S$-orbit), by proving topological versions of Dye's theorem \cite{GPS}. A key algebraic tool in this context is the dimension group, introduced by Elliott for the classification of AF-algebras \cite{E}: a dimension group is a direct limit of simplicial groups $\mathbb Z^{k_n}$ with positive connecting maps, abstractly characterized by Effros, Handelman and Shen \cite{EHS}, and visualized through Bratteli diagrams. The bridge to dynamical systems was built by Herman, Putnam and Skau \cite{HPS}, who proved that every essentially minimal system $(X,T)$ admits a Bratteli--Vershik realization, and that the dimension group of the associated diagram is isomorphic to the {\it ordered cohomology group} $K^0(X,T)=C(X,\mathbb Z)/\partial_TC(X,\mathbb Z)$ of the system,  where $\partial_T C(X,\mathbb Z)$ consists of the functions $f=h-h\circ T$ with $h$ continuous. This identifies the combinatorial data of the diagram with the dynamical $K$-theory. These studies remain active and are currently being studied to be extended to more general group actions \cite{GPS1}.

\medskip
Kakutani--Rokhlin (K-R) realizations of non-minimal zero-dimensional systems have since been studied systematically for aperiodic and densely aperiodic systems \cite{BDM, DK1, M, Tak2, Tak3}. Constructing K-R partitions for a general zero-dimensional system $(X,T)$ relies on a {\it quasi-section}: a closed subset $X_0\subseteq X$ such that every clopen neighbourhood of $X_0$ meets every orbit. Every point of a minimal system is a quasi-section, and quasi-sections exist in every zero-dimensional system \cite{M, Tak3}; in fact, Shimomura showed that every zero-dimensional system admits a {\it quasi-simple} quasi-section, that is, one meeting every minimal subset in exactly one point \cite{Tak3}. For a quasi-section $X_0$, the dimension group $K^0(V,E)$ of the associated Bratteli--Vershik realization is unital order isomorphic to the {\it relative dimension group}
$$K^0(X|X_0,T)=C(X,\mathbb Z)/\partial_T C(X|X_0,\mathbb Z),$$
where $\partial_T C(X|X_0,\mathbb Z)$ consists of the functions $f=h-h\circ T$ with $h$ continuous and constant on $X_0$ \cite{GW}. This relative structure is needed because, for non-transitive systems, $K^0(X,T)$ may fail to be partially ordered by the classes of non-negative functions, since there may be non-zero non-negative coboundaries.

\medskip
This paper addresses a natural generalization of a result of Glasner and Weiss \cite[Proposition 3.1]{GW}, which states that for a factor map $\pi:(X,T)\rightarrow (Y,S)$ between minimal Cantor systems, $K^0(Y,S)$ embeds in $K^0(X,T)$. For relative dimension groups one may expect a similar homomorphism
$$\pi^*: K^0(Y|Y_0,S) \rightarrow K^0(X|X_0,T),\qquad [f]\mapsto[f\circ\pi],$$
but the classical arguments rely on the Gottschalk--Hedlund theorem \cite{GotHed}, and on invariant measures, which are not adequate for non-minimal systems. We use instead the framework of \cite{GHY}, in which the factor map is represented as a limit of sliding block codes between Kakutani--Rokhlin realizations, so that $\pi^*$ is realized by a homomorphism $\iota$ between the corresponding direct limits (Proposition \ref{algebraic}).

\medskip
Our main result, Theorem \ref{injection}, states that if $X_0$ is a quasi-section of $(X,T)$ and $Y_0=\pi(X_0)$ is a quasi-simple quasi-section of $(Y,S)$, then $\pi^*$ is an injective unital order homomorphism. The key step (Lemma \ref{motor}) shows that if $f\circ\pi=h-h\circ T$ with $h$ vanishing on $X_0$, then $h$ vanishes on the whole preimage $\pi^{-1}(Y_0)$. The proof follows pairs of points in a common fibre into the minimal sets of $X$: along such pairs the difference of the values of $h$ is preserved, and quasi-simplicity forces these pairs to meet the fibre over the unique point of $Y_0$ in the corresponding minimal set of $Y$. If $h$ were not constant there, iterating this argument inside a minimal set would force $h$ to take values in an unbounded arithmetic progression, contradicting the boundedness of the continuous function $h$ (Lemma \ref{samemin}).

\medskip
We also characterize the factor maps for which $X_0=\pi^{-1}(Y_0)$ in terms of their premorphisms between Bratteli--Vershik realizations (Proposition \ref{equivalent}), and we show that in this case $\pi^*$ is a unital order embedding (Theorem \ref{equality}), with no quasi-simplicity assumption. In general, $\pi^*$ need not be an order embedding, even when $Y_0$ is quasi-simple, and the hypotheses of Theorem \ref{injection} cannot be dropped; we illustrate these phenomena with examples.

\medskip
The paper is organized as follows. Section 2 reviews some preliminaries on dynamical systems and ordered cohomology groups. In Section~\ref{sec:quasi}, we recall quasi-sections and Kakutani--Rokhlin partitions, which provide Bratteli--Vershik realizations of an arbitrary zero-dimensional system $(X,T)$, and we observe that a closed set is a quasi-section if and only if it meets every minimal subset. We then observe that the proof of \cite[Proposition~5.2]{GW}, stated there for minimal systems, uses minimality only to guarantee the existence of Kakutani--Rokhlin partitions; hence, once these are available with  quasi-sections, the same argument shows that for any quasi-section $X_0$ the relative dimension group $K^0(X|X_0,T)$ is unital order isomorphic to the dimension group of the associated Bratteli diagram (see also \cite{M} for the aperiodic case). As a consequence, the only non-negative conditional coboundary is zero, so no transitivity assumption is needed for $K^0(X|X_0,T)$ to be a partially ordered group. In Section~\ref{sec:factor}, we study how a factor map $\pi:(X,T)\to(Y,S)$ with $\pi(X_0)\subseteq Y_0$ acts on relative dimension groups via $\pi^*([g])=[g\circ\pi]$, and we identify $\pi^*$ with the homomorphism between the dimension groups induced by the premorphism of the associated Bratteli diagrams. We show that $\pi^*$ is a unital order embedding when $X_0=\pi^{-1}(Y_0)$ (Theorem~\ref{equality}), and that it is still an injective unital order-preserving homomorphism when $\pi(X_0)=Y_0$ and $Y_0$ is quasi-simple (Theorem~\ref{injection}); examples show that in the latter case $\pi^*$ need not be an order embedding and that quasi-simplicity cannot be dropped.

%%%%%%%%%%%%%%%%%%%%%%%%%%%%%%%%%%%%%%%%%%%%%%%%%%%%%%%
%%%%%%%%%%%%%%%%%%%%%%%%%%%%%%%%%%%%%%%%%%%%%%%%%%%%%%%%
\section{Preliminaries}
In this section, we briefly review essential definitions and prior results on zero-dimensional dynamical systems for our subsequent use. We refer the reader to the various resources  for the theory of Cantor systems and Kakutani-Rokhlin partitions  including \cite{DP, GH, GHY, HPS, P1, P2, Putnam, Tak2, Tak3} for the details and more information. 

 \subsection{zero-dimensional  dynamical systems.}\label{subsec_top}
A {\it zero-dimensional dynamical system} is a pair $(X, T)$, where $X$ is a non-empty compact totally disconnected metric space and $T$ is a homeomorphism of $X$. The topology of $X$ has a basis consisting of clopen (closed and open) sets. The orbit of a point $x\in X$, denoted by ${\mathcal O}(x)$, is the sequence $(T^nx)_{n\in\mathbb{Z}}$. If $X$ is finite, then $(X,T)$ is called a {\it trivial} dynamical system. For a point $x\in X$, the {\it $\omega$-limit set} is
$$\omega(x) = \{ y \in X :\ T^{n_k}x \to y \text{ for some sequence } n_k\to\infty \text{ in } \mathbb N \},$$
and the {\it $\alpha$-limit set} is
$$\alpha(x) = \{ y \in X :\ T^{-n_k}x \to y \text{ for some sequence } n_k\to\infty \text{ in } \mathbb N \}.$$

If $X$ is a Cantor space (that is, a non-empty, compact, metrizable, totally disconnected space with no isolated points), then the system is called a {\it Cantor system}. A system $(X,T)$ is {\it transitive} if it has at least one point with dense forward orbit, and it is {\it minimal} if all orbits are dense (equivalently, all forward orbits are dense). A non-trivial minimal zero-dimensional dynamical system is necessarily a Cantor system, since it has no isolated points. If the set of singleton fibres of $\pi$ is dense in $X$ then $\pi$ is called an {\it almost one-to-one} factoring. In minimal systems, if there is at least one point $y \in Y$ such that $\pi^{-1}(y)$ is a singleton, then $\pi$ is almost one-to-one. 

A subset $M\subseteq X$ is
called a \emph{minimal set} if $M$ is nonempty, closed and $T$-invariant
(i.e.\ $T(M)=M$), and no nonempty proper closed subset of $M$ is
$T$-invariant. In this case, $(M,T|_M)$ is called a \emph{minimal subsystem}
of $(X,T)$. For every $x\in X$, the set $\omega(x)$ contains a minimal set.

For $\varepsilon>0$, an {\it $\varepsilon$-chain} from $x$ to $y$ is a finite sequence $x=x_0,x_1,\ldots,x_n=y$ in $X$, with $n\geq 1$, such that $d(Tx_{i},x_{i+1})<\varepsilon$ for $0\leq i<n$. The system $(X,T)$ is {\it chain transitive} if for every $x,y\in X$ and every $\varepsilon>0$ there is an $\varepsilon$-chain from $x$ to $y$. Every transitive system is chain transitive, but not conversely; for instance, the one-point compactification $\mathbb Z\cup\{\infty\}$ of $\mathbb Z$ with the shift is chain transitive but not transitive.

Two topological dynamical systems $(X, T)$ and $(Y, S)$ are {\it semi-conjugate} if there exists a surjective continuous map $\pi: X\rightarrow Y$ such that $\pi\circ T=S\circ \pi$. In this case $(Y, S)$ is called a {\it factor} of $(X, T)$, $(X, T)$ is called an {\it extension} of $(Y, S)$, and $\pi$ is called a {\it factor map} or a {\it topological factoring}. Under a factor map, periodic points are mapped to periodic points, and every factor of a transitive (resp.\ minimal) system is transitive (resp.\ minimal).

 \subsection{Ordered Bratteli diagram}
 The  family of Bratteli-Vershik systems is correspondent  to the family of zero-dimensional dynamical systems \cite{HPS, Tak3}. Below are the essential materials  from this topic needed for the following sections.
 \begin{definition}\label{Bratteli}
\hspace{-0.2cm}
\begin{itemize}
\item A 
\emph{Bratteli diagram}
$B=(V, E)$ 
consists of an infinite sequence of finite, non-empty, pairwise disjoint sets 
$V_{0} = \{ v_{0} \}, V_{1}, V_{2}, \ldots$, 
called the vertices, another sequence of finite, non-empty, pairwise disjoint sets 
$E_{1}, E_{2}, \ldots$, 
called the 
\emph{edges}, 
and two maps 
$s : E_{n} \rightarrow V_{n-1}, r : E_{n} \rightarrow V_{n}$, 
for every 
$n \geq 1$, 
called the range and source maps, such that 
$r^{-1}\{ v \}$ 
is non-empty for all 
$v$ 
in 
$\cup_{n \geq 1} V_{n}$ 
and 
$s^{-1}\{ v \}$ 
is non-empty for all 
$v$ 
in 
$\cup_{n \geq 0} V_{n}$. For every $n\geq 1$ we have an {\it adjacency}  matrix $M_n$ of size $|V_n|\times |V_{n-1}|$ such that its entry $M_n^{ij}$ shows the number of edges between $v_i\in V_n$ and $v_j\in V_{n-1}$.

\item 
An 
\emph{ordered Bratteli diagram}
$B=(V, E, \leq)$
consists of a Bratteli diagram 
$(V, E)$ 
and a partial order
$\leq$
on 
$E$ 
such that two edges 
$e, e^{\prime}$ 
in 
$E$ 
are comparable if and only if 
$r(e) = r(e^{\prime})$. 
In such a diagram, we let 
$E_{\max}$ 
and 
$E_{\min}$ 
denote the set of maximal and minimal edges, respectively. If $0\leq m<n$ then $E_{m,n}$ denotes the set of all finite paths of edges from $V_m$ to $V_n$. Note that $E_{m,n}$ is a partially ordered set with respect to the reverse lexicographic ordering induced by $\leq$. 
\end{itemize}
\end{definition}

 \begin{definition}[\cite{BY, DK1}]\label{order}
Let  $B=(V,E,\leq)$ 
be an ordered Bratteli diagram and $X_B$ be the compact space of all (partially ordered) infinite paths of $B$.  Let $X_B^{\rm max}$ denote the set of all {\em infinite maximal paths} in $X_B$, i.e., the elements of $X_B$ each of whose component edges belonging to $E_{\max}$, and similarly let $X_B^{\rm min}$ denote the set of all {\it infinite minimal paths} in $X_B$. 

\begin{enumerate}
\item The Vershik map $T_B:X_B\setminus X_B^{\rm max} \rightarrow X_B\setminus X_B^{\rm min}$ is defined by
$$T_B(
e_1, e_2,
 \ldots, e_\ell,e_{\ell+1},\ldots)=(0,0,\ldots, e_{\ell}+1,e_{\ell+1},\ldots),$$
where $\ell$ is the first index such that $e_\ell$ is not the maximal edge in $r^{-1}(r(e_{\ell}))$,
$e_{\ell}+1$ is the successor of $e_\ell$
in  $r^{-1}(r(e_{\ell}))$,
 and  $0$ denotes the minimal edge in
  $r^{-1}(r(e_i))$  for every 
$i\geq 1$. 
\item $B$ is called \emph{perfect} if the Vershik map $T_B:X_B\setminus X_B^{\rm max}\rightarrow X_B\setminus X_B^{\rm min}$ can be extended to a continuous $T_B: X_B\rightarrow X_B$.

\item $B$  is called \emph{properly ordered} if it has a unique infinite minimal path $x_B^{\rm min}$. In this case, $T(X_B^{\rm max})=x_B^{\rm min}.$   
\end{enumerate}
\end{definition}

\subsection{Dimension group}\label{dimen}
A {\it unital\/} partially ordered group is a triple
$(G,\,G_+,\,u)$  where
\begin{itemize}
\item $G$ is an abelian group,
\item $G_+$ is a subset of $G$ such that
$$
G_+\cap (-G_+)=\{0\},\ \ \ G_++G_+\subseteq G_+,\ \ \ G_+-G_+=G.
$$
Note that $G_+$ with the above properties is called {\it the positive cone}.
\item $u$ is an element of $G_+$ such that for any $g\in G$ there exists
$n\in\mathbb Z^+$ such that $(nu-g)\in G_+$.
\end{itemize}
A  {\it dimension group} $G$ is a direct limit of finitely generated groups and positive homomorphisms:
$$\Z^{n_1}\stackrel{M_1} {\longrightarrow}\Z^{n_2}\stackrel{M_2}{\longrightarrow}\Z^{n_3}\stackrel{M_3}\longrightarrow\cdots\longrightarrow G=\varinjlim\Z^{(n_i)}$$
in which the groups $\Z^{n_i}$ are considered with their natural orderings and order units $(1,\ldots, 1)$. So a dimension group is unital and partially ordered. 
\medskip

For a partially ordered group $(G,G^+,u)$, a homomorphism $p:G\rightarrow \R$ is a  {\it state} (or {\it normalized  trace}) if $p(G^+)\geq 0$ and $p(u)=1$. The collection of all  states on $G$   is a compact convex subset of $\mathbb R^G$  denoted by $S(G,\,u)$. 
\medskip

\subsection{Ordered cohomology group}\label{cohom}
Let $(X,T)$ be a zero-dimensional system and let $C(X,\mathbb Z)$ be the group of continuous integer-valued functions on $X$. The subgroup of {\it coboundaries} is
\begin{equation}\label{cob}
\partial_T C(X,\mathbb Z)=\{f\in C(X,\mathbb Z):\ f=h-h\circ T \text{ for some } h\in C(X,\mathbb Z)\}.
\end{equation}
Two functions $f, g\in C(X,\mathbb Z)$ are called {\it cohomologous} if $f-g\in \partial_TC(X,\mathbb Z)$. Consider the quotient group
$$K^0(X, T) := C(X,\mathbb Z)/\partial_T C(X,\mathbb Z),$$
let
$$K^0_+(X,T)=\{[f]:\ f\in C(X,\mathbb Z),\ f\geq 0\}$$
be the set of classes of non-negative functions, and denote by $[1_X]$ the class of the constant function $1$.
%Let $\mathcal M_T(X)$ denote the set of $T$-invariant Borel probability measures on $X$, and $S(K^0(X,T),[1_X])$ the set of states of $K^0(X,T)$, that is, homomorphisms $\tau:K^0(X,T)\to\mathbb R$ with $\tau(K^0_+(X,T))\subseteq[0,\infty)$ and $\tau([1_X])=1$. There is a well-known bijection between $S(K^0(X,T),[1_X])$ and $\mathcal M_T(X)$, given by $\mu\mapsto\big([f]\mapsto\int f\,d\mu\big)$ \cite{BH, EHS}.
%\medskip
%A system $(X,T)$ is {\it transitive} if it has at least one point whose forward orbit is dense in $X$ (see \cite[Proposition 2.1.3]{DP}, where such systems are called recurrent).
If $f=h-h\circ T$ is a coboundary, then its Birkhoff sums telescope:
$$\sum_{i=0}^{n-1}f(T^ix)=h(x)-h(T^nx),$$
so they are bounded by $2\|h\|_\infty$, uniformly in $x$ and $n$. Suppose now that $(X,T)$ is transitive, and let $x_0\in X$ be a point whose forward orbit is dense. So, this forward orbit visits every non-empty open set infinitely often. Indeed, this is clear for an open set containing a non-isolated point; and if $p\in X$ is isolated, then, since $T$ is a homeomorphism, every $T^{-k}p$, $k\geq 0$, is isolated and hence lies on the forward orbit of $x_0$, which forces $x_0$ to be periodic, so that $X$ is a single periodic orbit. Now, if $f\geq 0$ and $f\neq 0$, then $f\geq 1$ on a non-empty clopen set, so $\sum_{i=0}^{n-1}f(T^ix_0)\to\infty$; thus $f$ is not a coboundary. Consequently, if $[f]\geq 0$ and $[-f]\geq 0$, say $f-g_1$ and $-f-g_2$ are coboundaries with $g_1,g_2\geq 0$, then $g_1+g_2$ is a non-negative coboundary, so $g_1=g_2=0$ and $[f]=0$. Hence
$$K^0_+(X,T)\cap\big(-K^0_+(X,T)\big)=\{0\},$$
and $(K^0(X,T), K^0_+(X,T), [1_X])$ is a partially ordered group with order unit, referred to as the {\it ordered cohomology group} of $(X,T)$ \cite{DP}. Boyle and Handelman showed that this group is partially ordered also for chain transitive systems \cite{BH}.

In general, non-zero non-negative coboundaries may exist, and then $K^0_+(X,T)$ is not a positive cone for $K^0(X,T)$. For examples, see Examples 4.2.4 and 4.3.2 in \cite{DP} and Example \ref{relative1} below.

\section{Quasi-sections and K-R partitions}\label{sec:quasi}
Using Krieger's marker techniques, the authors of \cite{DK1} proved that every zero-dimensional dynamical system is an inverse limit of subshifts. This can also be done with Kakutani--Rokhlin (K-R) partitions, as shown in \cite{Tak3}, following the methods of \cite{P1, Putnam}. This approach establishes, in particular, the existence of quasi-sections in every zero-dimensional dynamical system.

\begin{definition}[\cite{Tak3}]\label{basic}
Let $(X,T)$ be a zero-dimensional dynamical system.
\begin{itemize}
\item A closed subset $X_0\subseteq X$ is called a {\it quasi-section} if every clopen neighbourhood $U$ of $X_0$ is a complete $T$-section, i.e., $U$ meets every $T$-orbit of $X$ at least once; equivalently, $\bigcup_{n\in\mathbb Z} T^n(U)=X$.
\item A quasi-section $X_0$ is called {\it minimal} if it is minimal with respect to inclusion.
\item A quasi-section $X_0$ is called {\it quasi-simple} if for every minimal set $M$, $|X_0\cap M|=1$.
\end{itemize}
\end{definition}

The following lemma summarizes some direct consequences of the definition. We refer the reader to \cite{Tak3} for further results about quasi-sections.

\begin{lemma}[{\cite[Lemma 2.4]{Tak3}}]\label{minq}
Let $(X,T)$ be a zero-dimensional dynamical system and $X_0\subseteq X$ a closed set. The following are equivalent.
\begin{enumerate}
\item $X_0$ is a quasi-section for $(X,T)$.
\item For every $x\in X$, $\alpha(x)\cap X_0\neq\emptyset$.
\item For every $x\in X$, $\omega(x)\cap X_0\neq\emptyset$.
\item $X_0$ meets every minimal subset of $X$.
\end{enumerate}
\end{lemma}
\begin{proof}
The equivalence of (1), (2) and (3) is \cite[Lemma 2.4]{Tak3}. If (3) holds and $M$ is a minimal set, then for $x\in M$ we have $\omega(x)=M$, so $M\cap X_0\neq\emptyset$. Conversely, if (4) holds, then (3) follows, since every $\omega(x)$ contains a minimal set.
\end{proof}

We recall some of the arguments about the existence of Bratteli--Vershik (B-V) realizations for every zero-dimensional system.

\begin{theorem}[\cite{Tak3}]\label{zero1}
Every zero-dimensional dynamical system possesses non-trivial quasi-simple quasi-sections.
\end{theorem}

\begin{definition}\label{towers}
Let $(X,T)$ be a zero-dimensional dynamical system and $X_0\subseteq X$ a closed set. A {\it Kakutani--Rokhlin} (K-R) partition for $(X,T,X_0)$ is a partition
$$\mathcal P=\{Z(k,j):\ 1\leq k\leq K, \ 1\leq j\leq J(k)\}$$
of $X$ into clopen sets such that
\begin{enumerate}
\item $T(Z(k,j))=Z(k,j+1)$ for all $1\leq k\leq K$ and $1\leq j< J(k)$,
\item $T\big(\bigcup_{k=1}^K Z(k, J(k))\big)=\bigcup_{k=1}^K Z(k,1)$,
\item $X_0\subseteq \bigcup_{k=1}^K Z(k, 1)$.
\end{enumerate}
The set $\bigcup_{k=1}^K Z(k,J(k))$ is called the {\it top} of $\mathcal P$, and $\bigcup_{k=1}^K Z(k,1)$ its {\it base}, which we denote by $B(\mathcal P)$. For each fixed $k$, $1\leq k\leq K$, the set $\mathcal T=\{Z(k,j):\ 1\leq j\leq J(k)\}$ is called {\it a tower} of $\mathcal P$, with {\it base} $B(\mathcal T):=Z(k,1)$ and {\it height} $H(\mathcal T):=J(k)$.

\medskip
A {\it system of K-R partitions} $(\mathcal P_n)_{n\geq0}$ for $(X,T,X_0)$ is a sequence of K-R partitions such that, for every $n\geq 0$, $\mathcal P_{n+1}$ refines $\mathcal P_n$ and $B(\mathcal P_{n+1})\subseteq B(\mathcal P_n)$, and $\bigcup_{n\geq0}\mathcal P_n$ is a basis for the topology of $X$.
\end{definition}

\begin{theorem}[\cite{GHY}]\label{perfect}
Let $(X,T)$ be a zero-dimensional dynamical system and $X_0\subseteq X$ a closed non-empty set. The following statements are equivalent.
\begin{enumerate}
\item \label{per} There are a perfect ordered Bratteli diagram $B=(V,E,\leq)$ and a conjugacy $\gamma: (X,T)\rightarrow (X_B, T_B)$ such that $\gamma(X_0)=X_B^{\min}$.
\item \label{seq} There is a system of K-R partitions $(\mathcal P_n)_{n\geq0}$ for $(X, T,X_0)$ such that $\bigcap_{n\geq0} B(\mathcal P_n)=X_0$.
\item \label{wbase} $X_0$ is a quasi-section for $(X,T)$.
\end{enumerate}
\end{theorem}

By Theorem \ref{perfect}, for every zero-dimensional system $(X,T)$ and every quasi-section $X_0\subseteq X$ (which exists by Theorem \ref{zero1}), there exist systems of K-R partitions $(\mathcal P_n)_{n\geq 0}$ for $(X,T,X_0)$ leading to inverse limit realizations of $(X,T)$:
\[
\xymatrix{(\tilde{\mathcal P}_{0},\sigma)
&( \tilde{\mathcal P}_{1},\sigma)\ar[l]_{\alpha_{1}} &(\tilde{\mathcal P}_{2},\sigma)\ar[l]_{\alpha_{2}} &\cdots\ar[l]_{\alpha_3}&  (X,T,X_0),\ar[l]
}
\]
where $(\tilde{\mathcal P}_{n},\sigma)$ is the subshift associated with the K-R partition $\mathcal P_n$, which is a topological factor of $(X,T)$. The bases $B(\mathcal P_n)$ converge to $X_0$, and the set of infinite minimal paths of the associated Bratteli diagram corresponds bijectively to $X_0$ \cite[Section 5]{GHY}. It is well known that when $(X,T)$ is minimal, every singleton $\{x_0\}$, $x_0\in X$, is a quasi-section \cite{HPS}.

\medskip
Henceforth, we fix a quasi-section $X_0$ for each system $(X,T)$ and denote this structure by the triple $(X,T,X_0)$.

\subsection{Relative dimension group}\label{section2}
The ordered Bratteli-diagram $B=(V,E,\leq)$ in Theorem \ref{perfect} is called {\it a Bratteli-Vershik (B-V) realization } of $(X,T,X_0)$. Note that this realization leads to the existence of a dimension group, denoted  by $K^0(V,E)$,
\begin{equation}\label{dir}
\Z^{k_1}\stackrel{M_1} {\longrightarrow} \Z^{k_2}\stackrel{M_2} {\longrightarrow} \Z^{k_3}\stackrel{M_3} {\longrightarrow}\cdots\longrightarrow K^0(V,E)
\end{equation}
in which for every $n\geq 1$, $k_n=|V_n|$ and $M_n$ as in the Definition \ref{Bratteli}.  When $(X,T)$ is minimal, $K^0(V,E)$ is unital order isomorphic to the ordered cohomology group $K^0(X,T)$ \cite{DP, HPS}. While this isomorphism may fail in the non-minimal case, there exists a relative ordered cohomology group for the triple $(X,T,X_0)$ that remains unital order isomorphic to $K^0(V,E)$.

\medskip

 Let $(X,T)$ be a zero-dimensional dynamical system and $X_0\subseteq X$ a quasi-section. As in \cite{GW}, let $C(X|X_0,\mathbb Z)$ be the subgroup of $C(X,\mathbb Z)$ consisting of all functions that are constant on $X_0$. The subgroup of {\it conditional coboundaries} of $(X,T,X_0)$ is
$$\partial_T C(X|X_0,\mathbb Z)=\{f\in C(X,\mathbb Z):\ f=h-h\circ T \text{ for some } h\in C(X|X_0,\mathbb Z)\}.$$
Consider the quotient group
$$K^0(X|X_0,T):=C(X,\mathbb Z)/\partial_T C(X|X_0,\mathbb Z),$$
and let $K^0_+(X|X_0,T)=\{[f]:\ f\in C(X,\mathbb Z),\ f\geq 0\}$. Then $$(K^0(X|X_0,T),K^0_+(X|X_0,T),[1_X])$$ is a partially ordered group with order unit (see \cite{GPS, P1} for the original AF-algebraic arguments, and Theorem \ref{GW} and Corollary \ref{posneg} below). It is called the {\it relative dimension group} of $(X,T)$ corresponding to $X_0$ \cite{GW}.

In contrast with the ordered cohomology group of Subsection \ref{cohom}, no transitivity assumption is needed here: by Corollary \ref{posneg} (we will see below), the only non-negative function in $\partial_T C(X|X_0,\mathbb Z)$ is the zero function, so $K^0_+(X|X_0,T)\cap\big(-K^0_+(X|X_0,T)\big)=\{0\}$ for every zero-dimensional system and every quasi-section $X_0$. When $(X,T)$ is minimal, $X_0$ can be any single point of $X$ and in that case, $K^0(X|X_0,T)=K^0(X,T)$.
\medskip

The following theorem proves that $K^0(X|X_0,T)$ is unital order isomorphic to the dimension group $K^0(V,E)$, where $B=(V,E,\leq)$ is the B-V realization of $(X,T,X_0)$.
This theorem was first established for minimal systems in \cite[Proposition 5.2]{GW} and was confirmed to be hold for aperiodic systems in \cite{M}. By the definition of quasi-section and Theorem \ref{perfect}, the same proof as for \cite[Proposition 5.2]{GW} works for the following theorem. We present this proof, because it  is needed for the our discussions.
 
\begin{theorem}\label{GW}
Let $(X,T)$ be a zero-dimensional dynamical system and $X_0\subseteq X$ a quasi-section. Let $(Z_n)_{n\geq 0}$ be a descending sequence of clopen sets converging to $X_0$, with $Z_0=X$, and let $(R_n)_{n\geq 0}$ be a refining sequence of clopen partitions that forms a basis for the topology of $X$. Let $(\mathcal P_n)_{n\geq 0}$ be a refining sequence of K-R partitions which is compatible with $(R_n)_n$ and such that $Z_n$ serves as the base of $\mathcal P_n$. Denote by $k(n)$ the number of towers of $\mathcal P_n$ and by $H_i^{(n)}$ the height of its $i$-th tower, $1\leq i\leq k(n)$. Let $B = (V, E,\leq)$ be the corresponding ordered Bratteli diagram, so that $(X_B,T_B)$ is conjugate to $(X,T)$ with $X_{\min}$ corresponding to $X_0$, and let
\begin{equation}\label{direct}\xymatrix{
\mathbb{Z}\ar[r]^{M_0}&{\mathbb{Z}^{k(1)}} \ar[r]^{M_1}
&\mathbb{Z}^{k(2)}\ar[r]^{M_{2}}
&\cdots\ar[r]
&G=K^0(V,E)
}
\end{equation}
be the corresponding dimension group. Then $G$ is isomorphic (as a dimension group with order unit) to the relative dimension group $K^0(X|X_0, T)$.
In particular, the relative dimension group does not depend on the Bratteli-Vershik representation.
\end{theorem}

\begin{proof}
For $v\in\mathbb{Z}^{k(m)}$ we write $[v]_G$ for its image in $G$. The isomorphism is induced by a homomorphism $J\colon C(X,\mathbb{Z}) \rightarrow G$, defined as follows. Let $f\in C(X,\mathbb{Z})$. Since $f$ is locally constant and the partitions $R_n$ form a basis, there is $m$ such that $f$ is constant on each level of each tower of $\mathcal{P}_m$; let $f_m(i,j)$ denote the value of $f$ on the $j$-th level of the $i$-th tower of $\mathcal P_m$. Put
\begin{equation}\label{J}
\tilde{f}_m(i) = \sum_{j=0}^{H_i^{(m)}-1} f_m(i,j),\qquad 1\leq i\leq k(m),
\end{equation}
that is, the sum of the values of $f$ along the $i$-th tower of $\mathcal P_m$, and define
$$J(f)= [\tilde{f}_m]_G = [(\tilde{f}_m(1), \ldots, \tilde{f}_m(k(m)))]_G .$$
 
\emph{$J$ is well defined and onto.} Each tower of $\mathcal P_{m+1}$ is a concatenation of towers of $\mathcal P_m$, and the number of times the $i$-th tower of $\mathcal P_{m+1}$ passes through the $\ell$-th tower of $\mathcal P_m$ is the $(i,\ell)$ entry of $M_m$. Hence $\tilde f_{m+1}=M_m\tilde f_m$, so $J(f)$ does not depend on the choice of $m$, and $J$ is a homomorphism. If $v\in \mathbb{Z}^{k(m)}$ and $B_i^{(m)}$ denotes the base of the $i$-th tower of $\mathcal P_m$, then $f=\sum_i v_i 1_{B_i^{(m)}}$ satisfies $\tilde f_m=v$, so $J$ is onto.
 
To see that \emph{$\partial_TC(X|X_0,\mathbb{Z})\subseteq\ker(J)$}, let $f=h-h\circ T$ with $h\in C(X|X_0,\mathbb{Z})$. Subtracting from $h$ its constant value on $X_0$ does not change $h-h\circ T$, so we may assume $h(X_0)=\{0\}$. The set $U=\{h=0\}$ is a clopen set containing $X_0$. Since the $Z_n$ are compact, decreasing and $\bigcap_n Z_n=X_0\subseteq U$, we have $Z_m\subseteq U$ for all large $m$; choose such an $m$ for which $f$ is also constant on the levels of $\mathcal P_m$. If $x$ lies in the base of the $i$-th tower, then $x,\,T^{H_i^{(m)}}x\in Z_m$, and therefore
$$\tilde f_m(i)=\sum_{j=0}^{H_i^{(m)}-1}f(T^jx)=h(x)-h\big(T^{H_i^{(m)}}x\big)=0 .$$
Thus $J(f)=0$.
 
Now to see that \emph{$\ker(J)\subseteq\partial_TC(X|X_0,\mathbb{Z})$}, suppose $J(f)=0$. Then $[\tilde f_m]_G=0$, so $\tilde f_{m'}=M_{m'-1}\cdots M_m\tilde f_m=0$ for some $m'\geq m$. Replacing $m$ by $m'$, we may assume $\tilde f_m=0$, i.e.\ for every $x\in B(\mathcal P_m)$ lying in the $i$-th tower,
$$\sum_{j=0}^{H_i^{(m)}-1}f(T^jx)=0 .$$
So $h$ can be defined to be constant on the $\ell$-th level of the $i$-th tower of $\mathcal P_m$, $0\leq \ell\leq H_i^{(m)}-1$, with value
$$h\big|_{\text{level }\ell}=-\sum_{j=0}^{\ell-1}f_m(i,j)=\sum_{j=\ell}^{H_i^{(m)}-1}f_m(i,j).$$
Then $h$ is continuous and $h(B(\mathcal P_m))=0$; since $X_0\subseteq Z_m=B(\mathcal P_m)$, we get $h\in C(X|X_0,\mathbb{Z})$. If $x$ lies on level $\ell<H_i^{(m)}-1$, then $T(x)$ lies on level $\ell+1$ and $h(x)-h(T(x))=f_m(i,\ell)=f(x)$. If $x$ lies on the top level, then $T(x)\in B(\mathcal P_m)$, so $h(x)-h(T(x))=h(x)=f_m(i,H_i^{(m)}-1)=f(x)$. Hence $f=h-h\circ T$. Therefore, $f=h-h\circ T$.
 
Consequently $J$ induces a group isomorphism $K^0(X|X_0,T)\to G$. Moreover, $[f]\geq 0$ if and only if there exists some sufficiently large $m$ such that $\tilde{f}_m\geq 0$, so this isomorphism is an order isomorphism. Finally, since $Z_0=X$, the partition $\mathcal P_0$ consists of a single tower of height one, and $J(1_X)=[(H_1^{(m)},\ldots,H_{k(m)}^{(m)})]_G=[M_{m-1}\cdots M_0\, 1]_G$, which is the order unit of $G$. Thus the order unit $[1_X]$ is mapped to the order unit of $G$.
\end{proof}
 
\begin{corollary}\label{clopen}
Let $(X,T)$ be a zero-dimensional dynamical system and $X_0\subseteq X$ a quasi-section. Let $f\in \partial_TC(X|X_0,\mathbb{Z})$. There exists a clopen set $V$ containing $X_0$ such that for every $x\in V$ and every $n\in\mathcal N(x,V)$, $\sum_{i=0}^{n-1}f(T^ix)=0$. In particular, $f=h-h\circ T$ for some $h\in C(X|X_0,\mathbb{Z})$ which vanishes on $V$.
\end{corollary}
\begin{proof}
By Theorem \ref{GW}, $J(f)=0$, so, as in its proof, there is $m$ with $\tilde f_m=0$. Put $V=B(\mathcal P_m)=Z_m\supseteq X_0$. Let $x\in V$. Consider the elements of $\mathcal N(x,V)$ as a sequence:  $n_1<n_2<\cdots$ with $n_0=0$. Note that for each $k\geq 1$ and $n_{k-1}\leq s\leq n_k-1$, $T^s(x)$ passes through exactly one tower of $\mathcal P_m$ from its base to its top. So, for every $n_k\in \mathcal N(x,V)$ the sum $\sum_{i=0}^{n_k-1}f(T^ix)$ is a sum of $k$ entries of $\tilde f_m$, which are all zero. So the function $h$ constructed in the proof of Theorem \ref{GW} satisfies $f=h-h\circ T$ and vanishes on $V$.
\end{proof}
 
\begin{corollary}\label{posneg}
Let $(X,T)$ be a zero-dimensional dynamical system and $X_0 \subseteq X$ a quasi-section. If $f \in \partial_TC(X|X_0,\mathbb{Z})$ and $f \geq 0$, then $f \equiv 0$.
\end{corollary}
\begin{proof}
By Theorem \ref{GW}, $J(f)=0$, so, as in its proof, there is $m$ such that $f$ is constant on the levels of $\mathcal P_m$ and $\tilde f_m=0$. Each $\tilde f_m(i)$ is the sum of the values of $f$ on the levels of the $i$-th tower, and these values are non-negative; hence all of them vanish. Since the towers of $\mathcal P_m$ partition $X$, we conclude that $f\equiv 0$.
\end{proof}

%\subsection{Quasi-sections and conditional coboundaries}
%In light of Theorem \ref{GW}, we are led to investigate the connection between a system of K-R partitions for $(X,T,X_0)$ and the group of all coboundaries, extending the analysis beyond the conditional case. 
%
%
%

%\subsection{Invariant measures and quasi-sections}
%By Theorem \ref{GW}, if $X_0$ is a quasi-section for $X,T)$ there is a bijection between the sets of states on the two dimension groups $K^0(X|X_0,T)$ and $K^0(V,E)$, in which $B=(V,E,\leq)$ is the Bratteli-Vershik realization of $(X,T,X_0)$. Note that a state $\rho: K^0(X|X_0,T)\rightarrow \R$, is 
%%
%%%%%%%%%%%%%%%%%%%%%%%%%%%%%%%%%%%%%%%%%%%%%%%%%%%%%%%%%%%%%%%%%%%%%
%%%%%%%%%%%%%%%%%%%%%%%%%%%%%%%%%%%%%%%%%%%%%%%%%%%%%%%%%%%%%%%%%%%%%
%%%%%%%%%%%%%%%%%%%%%%%%%%%%%%%%%%%%%%%%%%%%%%%%%%%%%%%%%%%%%%%%%%%%%
%%%%%%%%%%%%%%%%%%%%%%%%%%%%%%%%%%%%%%%%%%%%%%%%%%%%%%%%%%%%%%%%%%%%%

\section{Topological factoring and relative dimension group}\label{sec:factor}
In this section, we consider two  zero-dimensional dynamical systems with a topological factoring between them $\pi:(X,T)\rightarrow (Y,S)$. We study the relation between the associated relative dimension groups of the two systems with respect to their non-trivial  quasi-sections $X_0\subseteq X$ and $Y_0\subseteq Y$, respectively when $\pi(X_0)\subseteq Y_0$.

% and $Y_0$ is  minimal.  We show that by Theorem \ref{inversefactor}, there is an order embedding  between the direct limit groups associated with the relative dimension groups $K^0(Y|Y_0,S)$ and $K^0(X|X_0,T)$ of the two systems.

%\subsection{Topological factoring and quasi-sections }
\medskip

We start by recalling the following theorem from \cite{GHY}.

 \begin{theorem}\label{inversefactor}
Consider zero-dimensional dynamical systems $(X,T,{ X_0})$ and $(Y, S, { Y_0})$,
where $X_0$ and $Y_0$ are quasi-sections. 
Then there exists a topological factoring
$ \pi:(X,T)\longrightarrow (Y,S)\ \ {\rm with} \ \  { \pi(X_0)\subseteq Y_0}$
if and only if  for every  K-R systems $(\mathcal P_n)_{n\geq 0}$ and $(\mathcal Q_n)_{n\geq 0}$ for $(X,T,X_0)$ and $(Y,S,Y_0)$, respectively, and the inverse limit systems associated to them,  there exists a strictly increasing sequence of natural numbers $\{n_i\}_{i\geq 0}$ and a sequence of sliding block codes $\pi_i:(\tilde{\mathcal P}_{n_i},\sigma)\rightarrow (\tilde{\mathcal Q}_i,\sigma)$, $i\geq 0$,  such that 
 all the following rectangles between the inverse limit sequences commute:
 \begin{equation}\label{diag}
\xymatrix{(\tilde{\mathcal{P}}_{0},\sigma)\ar[d]_{\pi_{0}}
&(\tilde{\mathcal{P}}_{n_1},\sigma)\ar[l]_{_{\gamma_{1}}}\ar[d]_{\pi_{1}} &(\tilde{\mathcal{P}}_{n_2},\sigma)\ar[l]_{\gamma_{2}}\ar[d]_{\pi_{2}} &\cdots\ar[l]_{\ \ \ \gamma_{3}}& (X,T,X_0)\ar[l]\ar[d]_{\pi}\ \ \  \\
(\tilde{\mathcal{Q}}_{0},\sigma)
&(\tilde{\mathcal{Q}}_1,\sigma)\ar[l]^{\beta_1} &(\tilde{\mathcal{Q}}_{2},\sigma)\ar[l]^{\beta_2}&\cdots\ar[l]^{\ \ \ \beta_3}&(Y,S, Y_0),\ar[l] \   
 }
\end{equation}
where 
$\gamma_i=\alpha_{n_i+1}\circ\alpha_{n_i+2}\circ\cdots\circ\alpha_{n_{i+1}}$.
\end{theorem}

The commutative diagram in Theorem \ref{inversefactor} induces the following relation between the Bratteli diagrams associated with the two systems. 
\begin{proposition}\cite{GHY}\label{localmor}
Let $(X,T)$ and $(Y,S)$ be zero-dimensional dynamical systems with quasi-sections $X_0$ and $Y_0$ respectively. Then there exists a toplogical factoring $\pi:X\rightarrow Y$ with $\pi(X_0)=Y_0$ if and only if  
for every B-V models $C=(W,E',\leq)$ and $B=(V,E,\leq)$ for $(X,T,X_0)$ and $(Y,S,\{y_0\})$ respectively,   there exists  an increasing sequence $\{n_i\}_{i\geq 0}$ of non-negative integers with $n_0=0$, and non-erasing letter-surjective morphisms $\eta_i:V_{n_i}^*\rightarrow W_i^*$  for every $i\geq 0$,  such that the following diagram commutes:
\begin{equation}\label{dig0}
\xymatrix{W_{n_i}^*\ar[d]_{\eta_{i}}
& \ \ W_{n_{i+1}}^*\ar[l]_{_{\sigma^C_{(n_i,n_{i+1}]}}}\ar[d]_{\eta_{i+1}} &\ \  \\
\ V_{i}^*
& \ \ V_{i+1}^*.\ar[l]^{\sigma^B_{i+1}} &\   
 }
\end{equation}
\end{proposition}

\begin{lemma}\label{quasi-section}
Let $\pi:(X,T)\rightarrow (Y,S)$ be a factor map between zero-dimensional dynamical systems.
\begin{enumerate}
\item If $X_0$ is a quasi-section for $(X,T)$, then $\pi(X_0)$ is a quasi-section for $(Y,S)$.
\item $Y_0\subseteq Y$ is a quasi-section for $(Y,S)$ if and only if $\pi^{-1}(Y_0)$ is a quasi-section for $(X,T)$.
\end{enumerate}
\end{lemma}
\begin{proof}
(1) The set $\pi(X_0)$ is compact, hence closed. Let $V$ be a clopen set containing $\pi(X_0)$. Then $U=\pi^{-1}(V)$ is a clopen set containing $X_0$, so $\bigcup_{n\in\mathbb Z}T^n(U)=X$. Since $\pi$ is onto, $\pi(U)=V$, and so
$$Y=\pi(X)=\bigcup_{n\in\mathbb Z}\pi(T^n(U))=\bigcup_{n\in\mathbb Z}S^n(\pi(U))=\bigcup_{n\in\mathbb Z}S^n(V).$$

(2) If $\pi^{-1}(Y_0)$ is a quasi-section, then, since $\pi$ is onto, so is $\pi(\pi^{-1}(Y_0))=Y_0$, by (1). Conversely, let $Y_0$ be a quasi-section; note that $\pi^{-1}(Y_0)$ is closed. Let $U$ be a clopen set containing $\pi^{-1}(Y_0)$, and write $Y_0=\bigcap_n V_n$ for a decreasing sequence of clopen sets $V_n$. Suppose that $\pi^{-1}(V_n)\not\subseteq U$ for every $n$, and pick $x_n\in\pi^{-1}(V_n)\setminus U$. Let $x$ be a limit point of $(x_n)_n$. For each $m$, we have $x_n\in\pi^{-1}(V_m)$ for all $n\geq m$, and this set is closed, so $\pi(x)\in\bigcap_m V_m=Y_0$; also $x\notin U$, since $X\setminus U$ is closed. This contradicts $\pi^{-1}(Y_0)\subseteq U$. Hence $\pi^{-1}(V)\subseteq U$ for some clopen set $V\supseteq Y_0$. Since $V$ is a complete $S$-section,
$$X=\pi^{-1}\Big(\bigcup_{n\in\mathbb Z}S^n(V)\Big)=\bigcup_{n\in\mathbb Z}T^n\big(\pi^{-1}(V)\big)\subseteq\bigcup_{n\in\mathbb Z}T^n(U),$$
so $U$ is a complete $T$-section. Hence $\pi^{-1}(Y_0)$ is a quasi-section.
\end{proof}

\begin{proposition}\label{equivalent}
Let $(X,T)$ and $(Y,S)$ be zero-dimensional dynamical systems with quasi-sections $X_0$ and $Y_0$, respectively. Suppose that $\pi:(X,T)\rightarrow (Y,S)$ is a factor map such that $\pi(X_0)\subseteq Y_0$. The following statements are equivalent.
\begin{enumerate}
\item $\pi^{-1}(Y_0)=X_0$.
\item For every B-V model $B$ for $(Y,S, Y_0)$ there exist a B-V model $C$ for $(X,T,X_0)$ and an ordered premorphism $f:B\rightarrow C$ such that $\pi=\alpha_f$ and, for every vertex $w$ of $C$, there exists a unique edge of $f$ with range $w$.
\item There exist B-V models $B$ and $C$ for $(Y,S,Y_0)$ and $(X,T,X_0)$, respectively, and an ordered premorphism $f:B\rightarrow C$ such that $\pi=\alpha_f$ and, for every vertex $w$ of $C$, there exists a unique edge of $f$ with range $w$.
\end{enumerate}
\end{proposition}
\begin{proof}
$(1)\Rightarrow (2)$. Suppose that $\pi^{-1}(Y_0)=X_0$, and let $B$ be a B-V model for $(Y,S,Y_0)$. Let $(\mathcal P_n)_{n\geq 0}$ be the associated nested sequence of K-R partitions of $(Y,S)$, with $\mathcal P_0=\{Y\}$, and let $Z_n$ be the base of $\mathcal P_n$. Then $Z_n\subseteq Z_{n-1}$ for $n\geq 1$ and $\bigcap_{n\geq 0}Z_n=Y_0$, so
$$\bigcap_{n\geq 0}\pi^{-1}(Z_n)=\pi^{-1}(Y_0)=X_0 .$$
Now for every $n\geq 0$, consider
$$\pi^{-1}(\mathcal P_n)=\{\pi^{-1}(A):\ A\in\mathcal P_n\}$$
which is a K-R partition of $(X,T)$ with base $\pi^{-1}(Z_n)$. Indeed, its towers are the sets $\pi^{-1}(\mathcal L)$, $\mathcal L\in\mathcal P_n$, with levels $T^j\big(\pi^{-1}(B(\mathcal L))\big)=\pi^{-1}\big(S^j(B(\mathcal L))\big)$.

Let $(\mathcal R_n)_{n\geq 0}$ be a refining sequence of clopen partitions of $X$ generating its topology. Put $\mathcal Q_0=\pi^{-1}(\mathcal P_0)=\{X\}$. Suppose that $n\geq 1$ and that K-R partitions $\mathcal Q_0,\ldots,\mathcal Q_{n-1}$ have been constructed so that, for $1\leq i\leq n-1$, $\mathcal Q_i$ refines $\mathcal Q_{i-1}$, $\pi^{-1}(\mathcal P_i)$ and $\mathcal R_i$, and the base of $\mathcal Q_i$ equals $\pi^{-1}(Z_i)$. Since $\pi^{-1}(Z_n)\subseteq\pi^{-1}(Z_{n-1})$, we can refine $\pi^{-1}(\mathcal P_n)$ to a K-R partition $\mathcal Q_n$ which refines $\mathcal Q_{n-1}$ and $\mathcal R_n$, by just partitioning the base $\pi^{-1}(Z_n)$ into finitely many clopen sets, each contained in a single level of $\mathcal Q_{n-1}$ and a single element of $\mathcal R_n$ along the corresponding return orbit, without changing the base. Then the base of $\mathcal Q_n$ is $\pi^{-1}(Z_n)$, the bases of $(\mathcal Q_n)_{n\geq 0}$ decrease to $X_0$, and $(\mathcal Q_n)_{n\geq 0}$ generates the topology of $X$. Let $C$ be the B-V model for $(X,T,X_0)$ associated with $(\mathcal Q_n)_{n\geq 0}$, and let $f:B\rightarrow C$ be the ordered premorphism induced by the inclusions of the towers of $\mathcal Q_n$ in those of $\pi^{-1}(\mathcal P_n)$. Therefore, $\pi=\alpha_f$.

Finally, since the base of $\mathcal Q_n$ coincides with the base $\pi^{-1}(Z_n)$ of $\pi^{-1}(\mathcal P_n)$, the first return time to $\pi^{-1}(Z_n)$ of each point in the base of a tower of $\mathcal Q_n$ is the height of that tower. Hence each tower of $\mathcal Q_n$ is contained in a single tower of $\pi^{-1}(\mathcal P_n)$. Therefore every vertex $w$ of $C$ is the range of exactly one edge of $f$.

\smallskip
$(2)\Rightarrow (3)$. This follows since $(Y,S,Y_0)$ admits a B-V model, by Theorem \ref{perfect} (see also \cite{GHY, Tak3}).

\smallskip
$(3)\Rightarrow (1)$. Let $B$, $C$ and $f$ be as in (3). Without loss of generality we assume that $X=X_C$, $Y=X_B$, $\pi=\alpha_f$, and that $T$ and $S$ are the Vershik maps. Then $X_0=X_C^{\min}$ and $Y_0=X_B^{\min}$. Since $\pi(X_0)\subseteq Y_0$, we have $X_0\subseteq \pi^{-1}(Y_0)$.

Conversely, let $x=(x_1,x_2,\ldots)\in\pi^{-1}(Y_0)$, so that $\pi(x)=(y_1,y_2,\ldots)\in X_B^{\min}$, and let $n\geq 1$. By order commutativity of $f$ and the isomorphism $E_n\circ F_n\simeq F_0\circ S_{0,f_n}$, the path $e_0y_1y_2\cdots y_{f_n}$ corresponds to the path $x_1x_2\cdots x_ne$, where $e_0$ and $e$ are edges of $f$. By (3), each of $e_0$ and $e$ is the unique edge of $f$ with its range, so both are minimal. Since $y_1,\ldots, y_{f_n}$ are minimal edges, the path $e_0y_1\cdots y_{f_n}$ is minimal; hence so is $x_1\cdots x_ne$, and therefore $x_1,\ldots,x_n$ are minimal edges. As $n$ is arbitrary, $x\in X_C^{\min}=X_0$. Thus $\pi^{-1}(Y_0)\subseteq X_0$, and so $\pi^{-1}(Y_0)=X_0$.
\end{proof}

%%%%%%%%%%%%%%%%%%%%%%%%%%%%%%%%%%%%%%%%%%%%%%%%%%%%%%%%%%%%%%%%%%%%%
%%%%%%%%%%%%%%%%%%%%%%%%%%%%%%%%%%%%%%%%%%%%%%%%%%%%%%%%%%%%%%%%%%%%%
%%%%%%%%%%%%%%%%%%%%%%%%%%%%%%%%%%%%%%%%%%%%%%%%%%%%%%%%%%%%%%%%%%%%%
%%%%%%%%%%%%%%%%%%%%%%%%%%%%%%%%%%%%%%%%%%%%%%%%%%%%%%%%%%%%%%%%%%%%%

Now we bring in a lemma, which is well known in K-theory, about the realization of homomorphisms between two dimension groups. It will be useful in the sequel.

\begin{lemma}\label{iota1}
Let $G$ and $H$ be dimension groups given by the direct limit representations in the following diagram, where all matrices $E_i$, $E'_i$ have non-negative entries. Suppose that for every $i\geq 1$, $F_i$ is a non-negative matrix, and that the diagram commutes, i.e.\ $F_1E'_1=E_1$ and $F_iE'_i=E_iF_{i-1}$ for every $i\geq 2$.
\begin{equation*}
\xymatrix{\mathbb Z\ar[d]_{\rm id}\ar[r]^{E'_1}
&\mathbb Z^{n_1}\ar[d]_{F_1}\ar[r]^{E'_2} &\mathbb Z^{n_2}\ar[d]_{F_2}\ar[r]^{E'_3}&\cdots\ar[r]& H\ar[d]_{\iota}  \\
\mathbb Z\ar[r]_{E_1}
&\mathbb Z^{m_1}\ar[r]_{E_2} &\mathbb Z^{m_2}\ar[r]_{E_3}&\cdots\ar[r]& G
 }
\end{equation*}
Let $\varphi_i:\mathbb Z^{n_i}\rightarrow H$ and $\psi_i:\mathbb Z^{m_i}\rightarrow G$ be the canonical homomorphisms. Define $\iota:H \rightarrow G$ by $\iota(z)=\psi_i(F_i(x))$, where $z=\varphi_i(x)$ for some $i\geq1$ and $x\in\mathbb Z^{n_i}$. Then $\iota$ is a well-defined, unital, order-preserving homomorphism.
\end{lemma}

\begin{proof}
For $i<j$ write $E'_{i,j}=E'_j\cdots E'_{i+1}$ and $E_{i,j}=E_j\cdots E_{i+1}$. By commutativity of the diagram, $F_jE'_{i,j}=E_{i,j}F_i$.

If $\varphi_i(x)=\varphi_j(y)$ with $i\leq j$, then there is $k\geq j$ such that $E'_{i,k}x=E'_{j,k}y$. Hence
$$\psi_i(F_ix)=\psi_k(E_{i,k}F_ix)=\psi_k(F_kE'_{i,k}x)=\psi_k(F_kE'_{j,k}y)=\psi_k(E_{j,k}F_jy)=\psi_j(F_jy),$$
so $\iota$ is well defined. It is a homomorphism, since each $F_i$ is linear and any two elements of $H$ can be represented at a common level.

The order units are $u_H=\varphi_1(E'_1(1))$ and $u_G=\psi_1(E_1(1))$, and $\iota(u_H)=\psi_1(F_1E'_1(1))=\psi_1(E_1(1))=u_G$, so $\iota$ is unital.

Finally,
$H^+ = \{ \varphi_i(x):\ i\geq 1,\ x \in (\mathbb Z^{n_i})^+ \},$
and since every $F_i$ is non-negative, $\iota(\varphi_i(x))=\psi_i(F_ix)\in G^+$ for $x\in(\mathbb Z^{n_i})^+$. Hence $\iota(H^+)\subseteq G^+$.
\end{proof}

Now consider two zero-dimensional systems $(X,T,X_0)$ and $(Y,S,Y_0)$ with B-V models $C=(W,E',\leq)$ and $B=(V,E,\leq)$, respectively, and suppose that $\pi:(X,T)\rightarrow (Y,S)$ is a factor map with $\pi(X_0)\subseteq Y_0$. In Proposition \ref{localmor}, we established that $\pi$ induces local morphisms between the levels of the two diagrams as shown in (\ref{dig0}). In \cite{AEG}, Elliott et al.\ proved that when the systems are minimal, there are an increasing sequence $(k_i)_{i\geq 0}$ and, for each $i\geq 0$, a morphism $F_i:V_i\rightarrow W_{k_i}$, in the opposite direction of $\eta_i: W_{k_i}\rightarrow V_i$, such that for every $i\geq 0$ the following diagram commutes:
\begin{equation}\label{dig3}
\xymatrix{V_{i}\ar[d]_{F_{i}}\ar[r]^{E_{i+1}}
& V_{i+1}\ar[d]^{F_{i+1}} \\
W_{k_i}\ar[r]_{E'_{(k_i,k_{i+1}]}}
& W_{k_{i+1}}
 }
\end{equation}
where $E'_{(k_i,k_{i+1}]}=E'_{k_{i+1}}\circ\cdots\circ E'_{k_i+1}$. The pair $f=\big((F_i)_{i\geq 0},(k_i)_{i\geq 0}\big)$ is called a {\it premorphism}. By Proposition \ref{localmor}, this construction extends to the non-minimal case as well. For $v\in V_i$ and $w\in W_{k_i}$, let $F_i(w,v)$ denote the number of edges of the morphism $F_i$ with source $v$ and range $w$. Thus $F_i$ determines the non-negative integer matrix $\big(F_i(w,v)\big)_{w\in W_{k_i},\,v\in V_i}$ of size $|W_{k_i}|\times|V_i|$. By abuse of notation, we denote this matrix again by $F_i$. Telescoping $C$ along $(k_i)_i$, Lemma \ref{iota1} yields a unital order-preserving homomorphism $\iota: K^0(V,E)\rightarrow K^0(W,E')$.

\begin{proposition}\label{algebraic}
Let $\pi:(X,T)\rightarrow (Y,S)$ be a factor map between zero-dimensional systems, and let $X_0$ and $Y_0$ be quasi-sections for $(X,T)$ and $(Y,S)$, respectively, such that $\pi(X_0)\subseteq Y_0$. Let $B$ and $C$ be B-V models for $(Y,S,Y_0)$ and $(X,T,X_0)$, respectively, and let $f:B\rightarrow C$ be the premorphism associated with $\pi$. Then the homomorphism
$$\pi^*: K^0(Y|Y_0,S)\rightarrow K^0(X|X_0, T),\qquad \pi^*([f])=[f\circ \pi],$$
is equivalent to the homomorphism $\iota$ of Lemma \ref{iota1}; that is, there are unital order isomorphisms $\alpha$ and $\beta$ such that the following diagram commutes:
\begin{equation}\label{dig7}
\xymatrix{K^0(Y|Y_0,S)\ar[d]_{\alpha}\ar[r]^{\pi^*}
& K^0(X|X_0,T)\ar[d]^{\beta} \\
K^0(V,E)\ar[r]_{\iota}
& K^0(W,E')
 }
\end{equation}
\end{proposition}

\begin{proof}
Let $(\mathcal P_i)_{i\geq 0}$ and $(\mathcal Q_k)_{k\geq 0}$ be the sequences of K-R partitions of $(Y,S)$ and $(X,T)$ associated with $B$ and $C$, and let $\alpha$ and $\beta$ be the isomorphisms $J$ of Theorem \ref{GW} for $(Y,S,Y_0)$ and $(X,T,X_0)$, respectively. By the construction of $f$ (see Theorem \ref{inversefactor} and \cite[Theorem 1.3]{GHY}), $\mathcal Q_{k_i}$ refines $\pi^{-1}(\mathcal P_i)$, and the edges of $F_i$ with range $w\in W_{k_i}$ correspond, in order, to the successive passages of the tower of $\mathcal Q_{k_i}$ associated with the vertex $w$ through the towers of $\pi^{-1}(\mathcal P_i)$. In particular, for every $v\in V_i$, the tower of $w$ passes exactly $F_i(w,v)$ times through the tower $\pi^{-1}(v)$ of $\pi^{-1}(\mathcal P_i)$.

Now, let $g\in C(Y,\mathbb Z)$, and choose $i$ such that $g$ is constant on the levels of $\mathcal P_i$. Then $g\circ\pi$ is constant on the levels of $\pi^{-1}(\mathcal P_i)$ and so on those of $\mathcal Q_{k_i}$. Therefore, (using the notations from the proof of Theorem~\ref{GW}) the sum of $g\circ\pi$ along a tower $w$ of $\mathcal Q_{k_i}$ is
$$\widetilde{(g\circ\pi)}_{k_i}(w)=\sum_{v}F_i(w,v)\,\tilde g_i(v),$$
that is, $\widetilde{(g\circ\pi)}_{k_i}=F_i\,\tilde g_i$. Therefore,
$$\beta(\pi^*[g])=\big[\widetilde{(g\circ\pi)}_{k_i}\big]=[F_i\tilde g_i]=\iota\big([\tilde g_i]\big)=\iota(\alpha[g]),$$
so the diagram commutes.
\end{proof}

 \medskip
 
\begin{theorem}\label{equality}
Let $\pi:(X,T)\rightarrow (Y,S)$ be a factor map between zero-dimensional systems. Suppose that $Y_0$ is a quasi-section for $(Y,S)$ and let $X_0=\pi^{-1}(Y_0)$. Then the homomorphism
$$\pi^*: K^0(Y|Y_0,S)\rightarrow K^0(X|X_0, T),\qquad \pi^*([g])=[g\circ \pi],$$
is a unital order embedding.
\end{theorem}

\begin{proof}
 By Lemma \ref{quasi-section}, $X_0$ is a quasi-section for $(X,T)$, and $\pi(X_0)=Y_0$. By Proposition \ref{equivalent}, for a B-V model $B=(V,E,\leq)$ of $(Y,S,Y_0)$ there exist a B-V model $C=(W,E',\leq)$ of $(X,T,X_0)$ and a premorphism $f=\big((F_i)_{i\geq 0},(k_i)_{i\geq 0}\big):B\rightarrow C$ with $\pi=\alpha_f$ such that every vertex of $C$ is the range of exactly one edge of $f$. By Proposition \ref{algebraic}, it suffices to prove the statement for the homomorphism $\iota:K^0(V,E)\rightarrow K^0(W,E')$ of Lemma \ref{iota1}. Let $\varphi_i:\mathbb Z^{|V_i|}\rightarrow K^0(V,E)$ and $\psi_k:\mathbb Z^{|W_k|}\rightarrow K^0(W,E')$ be the canonical homomorphisms, and for $i<j$ write $E_{(i,j]}=E_j\circ\cdots\circ E_{i+1}$. By the commutativity of (\ref{dig3}),
\begin{equation}\label{commute}
F_jE_{(i,j]}=E'_{(k_i,k_j]}F_i\qquad (i<j).
\end{equation}

Since every vertex $w\in W_{k_j}$ is the range of exactly one edge of $F_j$, every row of the matrix $F_j$ has exactly one entry equal to $1$ and all other entries equal to $0$. Also, every vertex $v\in V_j$ is the source of at least one edge of $F_j$ (because $\pi$ is onto), so every column of $F_j$ is non-zero as well. Consequently, $F_j$ has a non-negative integer left-inverse $L_j$. Hence, for every $y\in\mathbb Z^{|V_j|}$,
\begin{equation}\label{Fprop}
F_jy=0 \implies y=0, \qquad\text{and}\qquad F_jy\geq 0\implies y=L_jF_jy\geq 0 .
\end{equation}
So, $\iota$ is injective.  Because, if $x\in\mathbb Z^{|V_i|}$ with $\iota(\varphi_i(x))=\psi_{k_i}(F_ix)=0$ then, there is $j\geq i$ such that $E'_{(k_i,k_j]}F_ix=0$, and by (\ref{commute}), $F_j\big(E_{(i,j]}x\big)=0$. Now (\ref{Fprop}) implies that $E_{(i,j]}x=0$ and therefore, $\varphi_i(x)=\varphi_j(E_{(i,j]}x)=0$.

Since $F_i\geq 0$, the homomorphism $\iota$ preserves the ordering. Moreover, if $x\in\mathbb Z^{|V_i|}$ with $\iota(\varphi_i(x))=\psi_{k_i}(F_ix)\geq 0$ then, there is $j\geq i$ such that $E'_{(k_i,k_j]}F_ix\geq 0$, and by (\ref{commute}), $F_j\big(E_{(i,j]}x\big)\geq 0$. Then, by (\ref{Fprop}), $E_{(i,j]}x\geq 0$. So,   $\varphi_i(x)=\varphi_j(E_{(i,j]}x)\geq 0$.
Therefore $\iota(z)\geq 0$ if and only if $z\geq 0$ for every $z\in K^0(V,E)$ which means that $\iota$ is a unital order embedding.

\end{proof}
 
 Now we go through the case that we have a topological factoring between $(X,T)$ and $(Y,S)$ and for the quasi-section, $X_0\subseteq X$, $\pi(X_0)=Y_0$  but $X_0\neq \pi^{-1}(Y_0)$. 
 
 \medskip

 We first prove some  lemmas which will be used.
 
\begin{lemma}\label{shift}
Let $\pi:(X,T)\to(Y,S)$ be a factor map, $f\in C(Y,\mathbb Z)$ and
$h\in C(X,\mathbb Z)$ with $f\circ\pi=h-h\circ T$. Let $u,v\in X$ with
$\pi(u)=\pi(v)$.
\begin{enumerate}
\item[(a)] $h(T^nv)-h(T^nu)=h(v)-h(u)$ for all $n\in\mathbb Z$.
\item[(b)] If $(n_k)$ is a sequence of integers with $T^{n_k}u\to u'$ and
$T^{n_k}v\to v'$, then $\pi(u')=\pi(v')$ and $h(v')-h(u')=h(v)-h(u)$.
\end{enumerate}
\end{lemma}
\begin{proof}
\begin{itemize}
\item[(a)] Since $\pi(T^nu)=S^n\pi(u)=S^n\pi(v)=\pi(T^nv)$ for all $n$, we have
$$h(T^{n+1}v)-h(T^{n+1}u)=\big(h(T^nv)-f(\pi T^nv)\big)-\big(h(T^nu)-f(\pi T^nu)\big)
=h(T^nv)-h(T^nu).$$ Thus $n\mapsto h(T^nv)-h(T^nu)$ is constant on $\mathbb Z$, equal to its value
$h(v)-h(u)$ at $n=0$.

\item[(b)]$\pi(u')=\lim_k S^{n_k}\pi(u)=\lim_k S^{n_k}\pi(v)=\pi(v')$, and by (a) and
continuity of $h$,
$h(v')-h(u')=\lim_k\big(h(T^{n_k}v)-h(T^{n_k}u)\big)=h(v)-h(u)$.
\end{itemize}
\end{proof}

\begin{lemma}\label{samemin}
In the setting of Lemma~\ref{shift}, suppose that $M\subseteq X$  is a minimal set. If $a,b\in M$
with $\pi(a)=\pi(b)$, then $h(a)=h(b)$.
\end{lemma}
\begin{proof}
Let $D=h(b)-h(a)$, and put $u_0=a$, $u_1=b$. Suppose $u_j,u_{j+1}\in M$ have been
chosen with $\pi(u_j)=\pi(u_{j+1})$ and $h(u_{j+1})-h(u_j)=D$. Since $M$ is minimal,
the orbit of $u_j$ is dense in $M$, so there are integers $n_k$ with
$T^{n_k}u_j\to u_{j+1}$. Passing to a subsequence, $T^{n_k}u_{j+1}\to u_{j+2}$ for
some $u_{j+2}\in M$. Lemma~\ref{shift}(b), applied to
$u=u_j$, $v=u_{j+1}$, gives $\pi(u_{j+1})=\pi(u_{j+2})$ and
$h(u_{j+2})-h(u_{j+1})=h(u_{j+1})-h(u_j)=D$.
Inductively $h(u_j)=h(a)+jD$ for all $j\ge0$. But  $h$ is continuous and bounded, so $D=0$.
\end{proof}
 
\begin{lemma}\label{motor}
Let $\pi:(X,T)\to(Y,S)$ be a factor map of zero-dimensional systems and $X_0$  be a
quasi-section of $(X,T)$ with $Y_0=\pi(X_0)$ quasi-simple. Let $f\in C(Y,\mathbb Z)$
and $h\in C(X,\mathbb Z)$ with $h(X_0)=\{0\}$ and $f\circ\pi=h-h\circ T$. Then
$h(\pi^{-1}(Y_0))=\{0\}$.
\end{lemma}

\begin{proof}
Let $x_1\in\pi^{-1}(Y_0)$, choose $x_0\in X_0$ with $\pi(x_0)=\pi(x_1)$, and put
$D=h(x_1)-h(x_0)=h(x_1)$.

The set $\omega(x_1)$ contains a minimal set $M_b$; pick $b\in M_b$ and integers
$n_k$ with $T^{n_k}x_1\to b$. Passing to a subsequence, $T^{n_k}x_0\to a$ for some
$a\in X$. By Lemma~\ref{shift}(b), $\pi(a)=\pi(b)$ and $h(b)-h(a)=D$.

Next, $\omega(a)$ contains a minimal set $M_a$. Since $X_0$ is a quasi-section,
we may choose $a^*\in X_0\cap M_a$ and $b^*\in X_0\cap M_b$. As $a^*\in\omega(a)$,
there are integers $m_j$ with $T^{m_j}a\to a^*$; passing to a subsequence,
$T^{m_j}b\to b''$, and $b''\in M_b$ since $M_b$ is closed and invariant.
Again by Lemma~\ref{shift}(b), $\pi(b'')=\pi(a^*)$ and $h(b'')-h(a^*)=D$; since
$h(a^*)=0$, we get $h(b'')=D$.

Now $N:=\pi(M_a)$ is a minimal subset of $Y$, and $\pi(M_b)$ is a minimal set
containing $\pi(b'')=\pi(a^*)\in N$, so $\pi(M_b)=N$. Thus $\pi(a^*)$ and
$\pi(b^*)$ both lie in $Y_0\cap N$, and since $Y_0$ is quasi-simple,
$\pi(b^*)=\pi(a^*)=\pi(b'')$.

Finally, $b''$ and $b^*$ lie in the minimal set $M_b$ and have the same image
under $\pi$, so by Lemma~\ref{samemin}, $D=h(b'')=h(b^*)=0$.
\end{proof}

 \begin{theorem}\label{injection}
Let $\pi:(X,T)\to(Y,S)$ be a factor map of zero-dimensional systems, $X_0$ a
quasi-section of $(X,T)$, and suppose $Y_0=\pi(X_0)$ is a quasi-simple
quasi-section. Then $$\pi^*:K^0(Y|Y_0,S)\to K^0(X|X_0,T), \ \ \ \ \ \ \ \pi^*([f])=[f\circ\pi]$$
is an injective unital order preserving homomorphism.
\end{theorem}

\begin{proof}
Clearly $\pi^*$ is a homomorphism. Let $f\in C(Y,\mathbb Z)$ with
$g=f\circ\pi\in\partial_T C(X|X_0,\mathbb Z)$, say $g=h-h\circ T$ with $h$ constant
on $X_0$. Without loss of generality, we may assume that $h(X_0)=\{0\}$. We show
$f\in\partial_S C(Y|Y_0,\mathbb Z)$.
\medskip

By Lemma~\ref{motor}, $h(\pi^{-1}(Y_0))=\{0\}$. The set $W=h^{-1}(0)$ is clopen and
contains $\pi^{-1}(Y_0)$. Write $Y_0=\bigcap_n U_n$ for a decreasing sequence of
clopen sets. If $\pi^{-1}(U_n)\not\subseteq W$ for all $n$, pick
$x_n\in\pi^{-1}(U_n)\setminus W$; a limit point $x$ satisfies $\pi(x)\in Y_0$ and
$x\notin W$, a contradiction. So there is a clopen set $U\supseteq Y_0$ with
$\pi^{-1}(U)\subseteq W$, i.e.\ $h=0$ on $\pi^{-1}(U)$.

\medskip
Since $Y_0$ is a quasi-section, there is a K-R partition $\mathcal P$ of $(Y,S)$
with $Y_0\subseteq B(\mathcal P)\subseteq U$, refined so that $f$ is constant on
each of its levels. Let $\mathcal Q=\pi^{-1}(\mathcal P)$: its towers are the sets
$\pi^{-1}(\mathcal L)$, $\mathcal L\in\mathcal P$, with levels
$T^j\big(\pi^{-1}B(\mathcal L)\big)=\pi^{-1}\big(S^jB(\mathcal L)\big)$, and
$B(\mathcal Q)=\pi^{-1}(B(\mathcal P))\subseteq\pi^{-1}(U)$, so $h=0$ on
$B(\mathcal Q)$.
\medskip

Let $b,b'\in\pi^{-1}(B(\mathcal L))$ and $0\le j<H(\mathcal L)$. Since $f$ is
constant on the levels of $\mathcal L$, $g(T^tb)=f(S^t\pi b)=f(S^t\pi b')=g(T^tb')$
for $t<j$. Hence
$$h(T^jb)=h(b)-\sum_{t<j}g(T^tb)=0-\sum_{t<j}g(T^tb)=h(b')-\sum_{t<j}g(T^tb')=h(T^jb').$$

Now we can define $k\in C(Y,\mathbb Z)$ as it follows. On each  level $S^jB(\mathcal L)$, $j\geq 0$ of $\mathcal P$, $k$ can take the
common value of $h$ on $\pi^{-1}(S^jB(\mathcal L))$. Then $k$ is continuous (since the
levels are clopen), $k=0$ on $B(\mathcal P)\supseteq Y_0$, so
$k\in C(Y|Y_0,\mathbb Z)$, and $h=k\circ\pi$. Therefore
$$f\circ\pi=h-h\circ T=(k-k\circ S)\circ\pi,$$
and since $\pi$ is onto, $f=k-k\circ S\in\partial_S C(Y|Y_0,\mathbb Z)$. 
\medskip

Moreover, $\pi^*$ is order-preserving: if $[f]\geq 0$, then $f-g=k-k\circ S$ for some
$g\geq 0$ and $k\in C(Y|Y_0,\mathbb Z)$. Hence, 
$$f\circ\pi-g\circ\pi=k\circ\pi-(k\circ\pi)\circ T,$$
where $g\circ\pi\geq 0$ and $k\circ\pi\in C(X|X_0,\mathbb Z)$, so
$\pi^*[f]=[f\circ\pi]\geq 0$. Moreover, $1_Y\circ\pi=1_X$, so
$\pi^*[1_Y]=[1_X]$, and $\pi^*$ is unital.
\end{proof}
\medskip

\begin{remark}
In general $\pi^*$ need not be an order embedding as there may exists some $f\in C(Y,\Z)$ such that $[f]\not\geq 0$ while $[f\circ\pi]\geq 0$. 
\medskip

To see this, let $Y=\mathbb Z\cup\{\pm\infty\}$
with $T(k)=k+1,\ T(\pm\infty)=\pm\infty$, and let $Y_0=\{-\infty,0,1,+\infty\}$, and let $X$ consists of two copies
$\{A_n\}_{n\in\Z},\{B_n\}_{n\in\Z}$ of $\mathbb Z$ with common limits $r$ (backward) and $q$ (forward). Then $S$ is defined similarly on each copy. Now consider $\pi: X\rightarrow Y$ with $\pi(A_n)=\pi(B_n)=n$, $\pi(r)=-\infty$, $\pi(q)=+\infty$. Let $X_0=\{r,A_0,B_1,q\}$ and
$f=1_{\{-1\}}-1_{\{0\}}+1_{\{1\}}\in C(Y,\Z)$. We  have
$f\circ\pi-(1_{\{A_{-1}\}}+1_{\{B_1\}})=h-h\circ T$ with
$h=1_{\{A_1\}}-1_{\{B_0\}}\in C(X|X_0,\mathbb Z)$. So $\pi^*[f]=[f\circ\pi]\geq 0$, (since $1_{\{A_{-1}\}}+1_{\{B_1\}}\geq 0$)  whereas $[f]\not\geq 0$. Because, if $f-g=k-k\circ S$ with $g\geq 0$ and $k\in C(Y|Y_0,\mathbb Z)$, then, since $k$ is constant on $Y_0$, we have $f(0)-g(0)=k(0)-k(1)=0$, i.e.\ $g(0)=f(0)=-1<0$, a contradiction.
Note that here,  $X_0\neq \pi^{-1}(Y_0)=\{r,A_0, A_1, B_0, B_1, q\}$.
\end{remark}

\begin{example}\label{relative1}
An example of topological factoring between two non-minimal systems $\pi: (X,T) \rightarrow (Y,S)$, where the order embedding between the dimension groups $K^0(Y,S) \hookrightarrow K^0(X,T)$ fails, is provided by \cite[Examples 3.4, 3.5]{P2}. Here, we explore the order embedding 
$$\pi^*: K^0(Y|Y_0,S) \hookrightarrow K^0(X|X_0,T)$$ 
with respect to the non-trivial quasi-sections $X_0$ and $Y_0$ with $\pi(X_0) = Y_0$, where $Y_0$ is a quasi-simple (and minimal) quasi-section. Here $\pi^{-1}(Y_0)=X_0$.

\medskip
The two systems are defined by 
$$Y=\mathbb{Z}\cup\{\infty\},\ \ S(k)=k+1,\ S(\infty)=\infty,\ \  \\  
X=\mathbb{Z}\cup\{\infty,-\infty\},\ \ T(k)=k+1,\ T(\pm\infty)=\pm\infty.$$

The Bratteli Vershik realisation of $(X,T)$ is the stationary ordered diagram $B=(V,E,\leq)$ in which, at each level $n\geq 1$, the set of vertices is $V_n=\{v^{n}_1, v^n_2,v^n_3\}$ and the incidence matrix is $E_n=E=\begin{bmatrix}1&0&0\\ 1&1&1\\ 0&0&1\end{bmatrix}$. The partial ordering on the edges at each level is the left-to-right order. The two fixed points, positive and negative infinity, of $(X,T)$ correspond to the two vertical infinite minimum paths (which are also maximum paths) on the left and the right of the diagram, respectively; these are carried at each level by the vertices $v_1^n$ and $v_3^n$, respectively. Note that $X_0=X_{\text{min}}$ is the minimal quasi-section for $(X,T)$.

\medskip
Similarly, as $(Y,S)$ has one fixed point and  its Bratteli Vershik realization is the stationary ordered diagram $C=(V',E',\leq)$ in which, at each level $n\geq 1$, the set of vertices is $V_n=\{v'^{n}_1, v'^{n}_2\}$ and the incidence matrix is $E'_n=E'=\begin{bmatrix}1&0\\ 1&1\end{bmatrix}$. The partial ordering on the edges is the left-to-right ordering at every even level and the right-to-left ordering at every odd level. The unique fixed point, infinity, is the vertical infinite minimum path (which is also a maximum path) on the left that is carried at each level by the vertex $v'^{n}_1$. Note that $Y_0=Y_{\text{min}}$ is the minimal quasi-section for $(Y,S)$. 

\medskip
The relative dimension groups $G=K^0(X|X_0,T)$ and $H=K^0(Y|Y_0,S)=K^0(Y,S)$ are order isomorphic to $K^0(V,E)$ and $K^0(V',E')$, respectively, according to the direct limit systems and Theorem \ref{GW}. Therefore, for the system $(Y,S)$, all the computations are the same as in \cite[Example 3.5]{P2}, which shows that
$$H=\mathbb{Z}^2,\ \ H^+=\{(a,b):\ a>0 \ \text{or}\ a=0, \ b\geq 0\},\ \ u'=[1_Y]\sim (1,1).$$
However, there are some differences between $K^0(X,T)$ and $K^0(X|X_0,T)$. As described in \cite[Example 4.2.4]{DP}, there are non-zero non-negative functions in $C(X,\mathbb Z)$ which are coboundaries; for instance, $1_{\{-1\}}=h-h\circ T$ with $h=-1_{[0,+\infty]}$. Such an $h$ is not constant on $X_0$, since $h(-\infty)=0\neq -1=h(+\infty)$. In contrast, by Corollary \ref{posneg}, the only non-negative function in $\partial_T C(X|X_0,\mathbb{Z})$ is the zero function.
\medskip

To compute $G$, note that for every $a, b, c\in\mathbb{Z}$ and $n\geq 1$, we have $E^n(a,b,c)=(a, na+nc+b,c)$. In particular, the eventual kernel of the stationary direct limit system $\varinjlim\big(\mathbb{Z}^3\stackrel{E}{\longrightarrow}\mathbb{Z}^3\big)$ is trivial, so $G=\mathbb Z^3$, and the positive cone and the order unit are
$$G^+=\{(a,b,c):\ a\geq 0,\ c\geq 0,\ \text{and}\ (a+c>0\ \text{or}\ b\geq 0)\},\qquad u=[1_X]=(1,1,1).$$
The topological factoring $\pi:(X,T)\rightarrow (Y,S)$ is defined by
$$\pi(k)=k,\ \ \ \pi(\pm\infty)=\infty.$$ 
This leads to the following diagram and morphisms between the two direct limit systems:
 \begin{equation*}\label{dig2}
\xymatrix{\mathbb{Z}\ar[d]_{\text{id}}\ar[r]^{E'_1}
&\mathbb{Z}^{2}\ar[d]_{F_1}\ar[r]^{E'} &\mathbb{Z}^{2}\ar[d]_{F_2}\ar[r]^{E'}&\cdots\ar[r]& H\ar[d]_{\iota}  \\
\mathbb{Z}\ar[r]_{E_1}
&\mathbb{Z}^{3}\ar[r]_{E}&\mathbb{Z}^{3}\ar[r]_{E}&\cdots\ar[r]& G,
 }
\end{equation*}
in which $E'_1=\begin{bmatrix} 1\\1\end{bmatrix}$, $E_1=\begin{bmatrix} 1\\ 1\\ 1\end{bmatrix}$, $E$ and $E'$ are as defined above, and for every $n\geq 1$, $F_n=\begin{bmatrix}1&0\\ n-1&1\\ 1&0\end{bmatrix}$. It is clear that for every $n\geq 1$, the kernel of $F_n$ is $\{(0,0)\}$, and so $\iota$ is injective. Note that $\iota$ is unital and $\iota(H^+)\subset G^+$ because for $(a,b)\in H^+$, there exists some $n\geq 1$ such that $F_n(a,b)=(a, (n-1)a+b, a)\geq 0$. 
\end{example}

\begin{example}\label{relative2}
For examples of having an injective $\pi^*$ in the case where $\pi(X_0)=Y_0$ but $\pi^{-1}(Y_0)\neq X_0$, consider any topological factoring $\pi:(X,T)\to(Y,S)$ between two minimal Cantor systems which is not almost one-to-one, for instance, \cite[Example 3.4]{GH}. Fix $y_0\in Y$ and $x_0\in\pi^{-1}(y_0)$, and choose Bratteli--Vershik realizations of $(Y,S)$ and $(X,T)$ whose unique infinite minimal paths correspond to $y_0$ and $x_0$, respectively. By minimality, $Y_0=\{y_0\}$ and $X_0=\{x_0\}$ are minimal quasi-sections, $Y_0$ is quasi-simple, and $\pi(X_0)=Y_0$. On the other hand, $\pi^{-1}(Y_0)\neq X_0$, since a factor map between minimal systems which is not almost one-to-one has no singleton fibres.

Since $X_0$ and $Y_0$ are singletons, $K^0(X|X_0,T) = K^0(X,T)$ and $K^0(Y|Y_0,S) = K^0(Y,S)$. Hence, by Theorem \ref{injection} (see also \cite[Proposition 3.1]{GW}), $\pi^*$ is an injective unital order embedding.
\end{example} 

\begin{remark}
In Theorem \ref{injection}, neither the quasi-simplicity of $Y_0$ nor the condition $\pi(X_0)=Y_0$ can be dropped: without them, the homomorphism $\iota$, and hence $\pi^*$, need not be injective. In general, since every closed superset of a quasi-section is again a quasi-section, one could consider quasi-sections with $\pi(X_0)\subsetneq Y_0$ and no constraint on $Y_0\setminus\pi(X_0)$.  In that generality, comparing $K^0(Y|Y_0,S)$ and $K^0(X|X_0,T)$ is of little interest.

 For a simple counterexample, let $X=\{x_1, x_2, x_3, x_4\}$ with $T(x_1)=x_2$, $T(x_2)=x_1$, $T(x_3)=x_4$, $T(x_4)=x_3$, and let $Y=\{y_1, y_2\}$ with $S(y_1)=y_2$, $S(y_2)=y_1$. The map $\pi:X\rightarrow Y$ given by $\pi(x_1)=\pi(x_3)=y_1$ and $\pi(x_2)=\pi(x_4)=y_2$ is a factor map. Let $Y_0=Y$, which is a quasi-section but not quasi-simple, and let $X_0$ be  $\{x_1, x_3\}$. Then, $\pi(X_0)=\{y_1\}\subsetneq Y_0$. So the quasi-simplicity fails.
%In both cases, the direct limit systems of Proposition \ref{algebraic} are obtained with $n_i=m_i=2$ for every $i\geq 1$, and
%$$E_1=\begin{bmatrix}2\\2\end{bmatrix},\quad E_i={\rm id}_{2\times 2}\ (i\geq 2),\qquad E'_1=\begin{bmatrix}1\\1\end{bmatrix},\quad E'_i={\rm id}_{2\times 2}\ (i\geq 2),$$
%where $E_i$ and $E'_i$ correspond to $(X,T)$ and $(Y,S)$, respectively. The connecting maps are $F_i=\begin{bmatrix}1&1\\1&1\end{bmatrix}$ for every $i\geq 1$, and $\ker F_i=\{(a,-a):\ a\in\mathbb Z\}\cong\mathbb{Z}$. Since all $E'_i$ are injective, $\iota$ is not injective.
%
Now, let $f=1_{\{y_1\}}-1_{\{y_2\}}$. Since $Y_0=Y$, every $k\in C(Y|Y_0,\mathbb Z)$ is constant, so $\partial_SC(Y|Y_0,\mathbb Z)=\{0\}$ and $[f]\neq 0$. On the other hand, $f\circ\pi=h-h\circ T$, where $h=1_{\{x_1,x_3\}}$. So,  $h$ is constant on $X_0$. Hence $\pi^*[f]=0$.
\end{remark}

{\bf Acknowledgements.}  The author is grateful to Nasser Golestani for some helpful discussions during the early stages of this work, and for his questions, which helped clarify some of the results. Lemma~\ref{iota1} arose from these discussions. The author also thanks Thierry Giordano for helpful comments on a version of the introduction.

%%%%%%%%%%%%%%%%%%%%%%%%%%%%%%%%%%%%%%%%%%%%%%%%%

\end{document}